\documentclass[11pt]{amsproc}
\usepackage{amssymb,float,amsmath,times}

\theoremstyle{plain}
\newtheorem{thm}{Theorem}[section]
\newtheorem{lem}[thm]{Lemma}

\newcommand{\PG}{\mathsf{PG}}
\newcommand{\GF}{\mathsf{GF}}
\newcommand{\Z}{\mathbb{Z}}
\newcommand{\f}{\mathcal F}
\newcommand{\ff}{\mathcal F^*}
\newcommand{\lcs}[2]{\gamma_{#2}(#1)}
\newcommand{\F}{\mathbb F}

\renewcommand{\le}{\leqslant}
\renewcommand{\ge}{\geqslant}

\begin{document}

\title{Elation generalised quadrangles of order $(s,p)$, where $p$ is prime}
 
\author[Bamberg]{John Bamberg}
\address{ %
Department of Pure Mathematics,
Ghent University,
Galglaan 2, B-9000 Ghent, Belgium.}
\email{bamberg@cage.ugent.be}

\author[Penttila]{Tim Penttila}
\address{ %
Department of Mathematics,
Colorado State University,
Fort Collins, CO 80523,
USA.}
\email{penttila@math.colostate.edu}

\author[Schneider]{Csaba Schneider}
\address{ %
  Informatics Research Laboratory, Computer and Automation Research
  Institute, 1518 Budapest Pf. 63, Hungary.}
\email{csaba.schneider@sztaki.hu}



\begin{abstract}
  We show that an elation generalised quadrangle which has
  $p+1$ lines on each point, for some prime $p$, is classical or arises
  from a flock of a quadratic cone (i.e., is a \textit{flock quadrangle}). 
\end{abstract}

\maketitle

\section{Introduction}

A \textit{generalised quadrangle} is an incidence structure of points and
lines such that if $P$ is a point and $\ell$ is a line not incident with
$P$, then there is a unique line through $P$ which meets $\ell$ in a
point. From this property, one can see that there are constants $s$ and $t$
such that each line is incident with $t+1$ points, and each point is
incident with $s+1$ lines. Such a generalised quadrangle is said to have
\emph{order} $(s,t)$, and hence its point-line dual is a generalised
quadrangle of order $(t,s)$. Of the known generalised quadrangles, most
admit a group of elations (see Section \ref{background} for a definition)
and are called \emph{elation generalised quadrangles}. In this paper, we
will be interested in elation generalised quadrangles where the parameter
$t$ is prime.

If $\mathcal{S}$ is an elation generalised quadrangle of order $(p,t)$, for
some prime $p$, then the elation group $G$ is a $p$-group (see \cite[Lemma
6]{Fro88}, and note that $s$ and $t$ are interchanged!).  In this
situation, we have by a deep result of Bloemen, Thas, and Van Maldeghem
\cite{BTvM}, that $\mathcal{S}$ is isomorphic to one of the classical
generalised quadrangles $\mathsf{W}(p)$, $\mathsf{Q}(4,p)$, or
$\mathsf{Q}(5,p)$. The same is not true if we interchange points and
lines. Suppose that $\mathcal{S}$ is an elation generalised quadrangle of
order $(s,p)$ (where $p$ is a prime).  Again, by a result of Frohardt
\cite[Lemma 6]{Fro88}, we have that the elation group is a $p$-group,
however, there exist candidates for $\mathcal{S}$ which are not classical
but are known as \textit{flock quadrangles}. These elation generalised
quadrangles are obtained from a flock of $\PG(3,p)$ (a partition of the
points of a quadratic cone of $\PG(3,p)$, minus its vertex, into conics)
and they have order $(p^2,p)$.  Such a quadrangle is classical if and only
if the flock is linear; and there do exist non-linear flocks for $p$ a prime at least
$3$ \cite[\S 10.6]{FGQ}. In this paper, we prove a result that is complementary
to that of Bloemen, Thas, and Van Maldeghem:

\begin{thm}\label{Theorem1}
  If $p$ is a prime, then an elation generalised quadrangle of order
  $(s,p)$ is classical or a flock quadrangle.
\end{thm}

Note that the above theorem does not hold when $p$ is replaced by a prime power
since the duals of the Tits quadrangles $T_3(O)$ arising from the Tits ovoids are 
elation generalised quadrangles of order $(q^2,q)$ (for $q=2^h$ and $h$ an odd number at least $3$)
that are not flock quadrangles, and the Roman elation generalised quadrangles of Payne are 
of order $(q^2,q)$ (for $q=3^h$, $h>2$) but are not flock quadrangles.

The proof of Theorem~\ref{Theorem1} relies on the following result
concerning Kantor families for groups of order $p^5$ (see
Section~\ref{background} for a definition of Kantor families).

\begin{thm}\label{Csaba}
  If $p$ is an odd prime and $G$ is a finite $p$-group of order $p^5$ that
  admits a Kantor family of order $(p^2,p)$, then $G$ is an extraspecial
  group of exponent $p$.
\end{thm}

In Sections \ref{background}--\ref{flocksec}, we briefly revise the basic
background theory and definitions needed for this paper.  Kantor families
for groups of order $p^5$ are then investigated in Section~\ref{sectp5},
and Theorem~\ref{Csaba} is proved in Section~\ref{proofthcsaba}.  Finally
in Section \ref{proofs} we prove Theorem \ref{Theorem1}.

Though our group theoretic notation is standard, we briefly review it for
the sake of a reader whose interest lies more in geometry than in group
theory.  If $a$ is a group element of order $p$ and $\alpha\in\F_p$ then,
identifying $\alpha$ with an element in $\{0,\ldots,p-1\}$, we may write
$a^\alpha$.  If $a$ and $b$ are group elements, then we define their {\em
  commutator} as $[a,b]=a^{-1}b^{-1}ab$. The properties of group
commutators that we need in this paper are listed, for instance,
in~\cite[\S~5.1.5]{robinson}.  The {\em centre} of a group $G$ consists of
those elements $z\in G$ that satisfy $[g,z]=1$ for all $g\in G$.  If $H,\
K$ are subgroups of a group $G$, then the {\em commutator subgroup} $[H,K]$
is generated by all commutators $[a,b]$ where $a\in H$ and $b\in K$. The
{\em derived subgroup} $G'$ of $G$ is defined as $[G,G]$.  The symbol $\lcs
Gi$ denotes the $i$-th term of the {\em lower central series} of $G$; that
is $\lcs G1=G$, $\lcs G2=G'$, and, for $i\ge 3$, $\lcs G{i+1}=[\lcs Gi,G]$.
The {\em nilpotency class} of a $p$-group is the smallest $c$ such that
$\lcs G{c+1}=1$.  The {\em Frattini subgroup} $\Phi(G)$ of a finite group $G$
is the intersection of all the maximal subgroups. If $G$ is a finite
$p$-group, then $\Phi(G)=G'G^p$ and $\log_p|G:\Phi(G)|$ is the size of a
minimal set of generators for $G$.  The basic properties of the Frattini
subgroup of a $p$-group can be found, for instance,
in~\cite[\S~5.3]{robinson}.  The {\em exponent} of a finite group $G$ is
the smallest positive $n$ such that $g^n=1$ for all $g\in G$.

%
%

\section{Generalised Quadrangles and Kantor Families}\label{background}

\subsection*{The basics}

A (finite) generalised quadrangle is an incidence structure of points
$\mathcal{P}$, lines $\mathcal{L}$, together with a symmetric point-line
incidence relation satisfying the following axioms:
\begin{enumerate}
\item[(i)] Each point lies on $t+1$ lines ($t\ge 1$) and two distinct
  points are incident with at most one line.
\item[(ii)] Each line lies on $s+1$ points ($s\ge 1$) and two distinct
  lines are incident with at most one point.
\item[(iii)] If $P$ is a point and $\ell$ is a line not incident with $P$,
  then there is a unique point on $\ell$ collinear with $P$.
\end{enumerate}
We say that our generalised quadrangle has order $(s,t)$ (or order $s$ if
$s=t$), and the point-line dual of a generalised quadrangle of order
$(s,t)$ is again a generalised quadrangle but of order $(t,s)$.  Higman's
inequality states that the parameters $s$ and $t$ bound one another; that
is, $t\le s^2$ and dually, $s\le t^2$.  A collineation $\theta$ of
$\mathcal{S}$ is an \textit{elation} about the point $P$ if it is either
the identity collineation, or it fixes each line incident with $P$ and
fixes no point not collinear with $P$. If there is a group $G$ of elations
of $\mathcal{S}$ about the point $P$ such that $G$ acts regularly on the
points not collinear with $P$, then we say that $\mathcal{S}$ is an
\textit{elation generalised quadrangle} with elation group $G$ and
\textit{base point} $P$. Necessarily, $G$ has order $s^2t$.

The \emph{classical generalised quadrangles} $\mathsf{W}(q)$,
$\mathsf{Q}(4,q)$, $\mathsf{H}(3,q^2)$, $\mathsf{Q}(5,q)$ and
$\mathsf{H}(4,q^2)$, are elation generalised quadrangles and arise as polar
spaces of rank $2$.  The first of these is the incidence structure of all
totally isotropic points and totally isotropic lines with respect to a null
polarity in $\PG(3,q)$, and is a generalised quadrangle of order $q$. The
point-line dual of $\mathsf{W}(q)$ is $\mathsf{Q}(4,q)$, the parabolic
quadric of $\PG(4,q)$, and is therefore a generalised quadrangle of order
$q$ (see \cite[3.2.1]{FGQ}). The incidence structure of all points and
lines of a non-singular Hermitian variety in $\PG(3,q^2)$, which forms the
generalised quadrangle $\mathsf{H}(3,q^2)$ of order $(q^2,q)$, has as its
point-line dual the elliptic quadric $\mathsf{Q}(5,q)$ in $\PG(5,q)$, which
is a generalised quadrangle of order $(q,q^2)$ (see \cite[3.2.3]{FGQ}). The
remaining classical generalised quadrangle, $\mathsf{H}(4,q^2)$, is the
incidence structure of all points and lines of a non-singular Hermitian
variety in $\PG(4,q^2)$, and is of order $(q^2,q^3)$ (see
\cite[3.1.1]{FGQ}).

\subsection*{Kantor families}

Now standard in the theory of elation generalised quadrangles are the
equivalent objects known commonly as \textit{$4$-gonal families} or
\textit{Kantor families} (after their inventor).  Let $G$ be a group of
order $s^2t$ and suppose there exist two families of subgroups
$\mathcal{F}=\{A_0,\ldots,A_t\}$ and $\mathcal{F}^*=\{A_0^*,\ldots,A_t^*\}$
of $G$ such that
\begin{enumerate}
\item[(a)] every element of $\mathcal{F}$ has order $s$ and every element
  of $\mathcal{F}^*$ has order $st$;
\item[(b)] $A_i\le A_i^*$ for all $i$;
\item[(c)] $A_i\cap A_j^*=1$ for $i\ne j$ (the ``tangency condition'');
\item[(d)] $A_iA_j\cap A_k=1$ for distinct $i$,$j$,$k$ (the ``triple condition'').
\end{enumerate}
Then the triple $(G,\mathcal{F},\mathcal{F}^*)$ is called a \textit{Kantor
  family}, but we will also say that $(\mathcal{F},\mathcal{F}^*)$ is a
\textit{Kantor family for $G$}. The pair $(s,t)$ is said to be the {\em
  order} of $(\f,\ff)$. From a Kantor family as described above, we can
define a point-line incidence structure as follows.

\begin{table}[H]
\begin{center}
\begin{tabular}{l|l}
\hline
Points & Lines \\
\hline
elements $g$ of $G$&
the right cosets $A_ig$\\
right cosets $A_i^*g$&
symbols $[A_i]$\\
a symbol $\infty$.&\\
\hline
\end{tabular}
\end{center}
\caption{The points and lines of the elation generalised quadrangle arising
  from a Kantor family (n.b., $A_i\in \mathcal{F}$, $A_i^*\in\mathcal{F}^*$, $g\in G$).}
\end{table}

Incidence comes in four flavours (points on the left, lines on the right):
\begin{center}
\begin{tabular}{rcl}
$g$&$\sim $&$A_ig$\\
$A_i^*g $&$\sim$&$[A_i]$\\
$A_i^*g$&$\sim$&$A_ih$, where $A_ih\subseteq A_i^*g$\\
$\infty$&$\sim$&$[A_i]$
\end{tabular}
\end{center}

It turns out that this incidence structure is an elation generalised
quadrangle of order $(s,t)$ with base point $\infty$ and elation group
$G$. Remarkably, all elation generalised quadrangles arise this way
\cite[\S 8.2]{FGQ}, and we obtain a so-called \textit{translation
  generalised quadrangle} when $G$ is abelian~\cite[8.2.3]{FGQ}.

%
%

\section{Flock Quadrangles and Special Groups}\label{flocksec}

\subsection*{Flock generalised quadrangles}\label{flocks}

A \textit{flock} of the quadratic cone $\mathcal{C}$ with vertex $v$ in
$\PG(3,q)$ is a partition of the points of $\mathcal{C}\backslash\{v\}$
into conics. Thas \cite{Thas87} showed that a flock gives rise to an
elation generalised quadrangle of order $(q^2,q)$, which we call a
\textit{flock quadrangle}.  The flocks of $\mathsf{PG}(3,q)$ have been
classified by Law and Penttila \cite{LawPenttila03} for $q$ at most $29$.
A \textit{BLT-set of lines} of $\mathsf{W}(q)$ is a set $\mathcal{L}$ of
$q+1$ lines of $\mathsf{W}(q)$ such that no line of $\mathsf{W}(q)$ is
concurrent with more than two lines of $\mathcal{L}$. For $q$ odd, Knarr
\cite{Knarr92} gave a direct geometric construction of an elation
generalised quadrangle from a BLT-set of lines of $\mathsf{W}(q)$.  The
ingredients of the Knarr construction are as follows:
\begin{itemize}
\item a symplectic polarity $\rho$ of $\mathsf{PG}(5,q)$;
\item a point $P$ of $\mathsf{PG}(5,q)$;
\item a 3-space inducing a $\mathsf{W}(q)$ contained in $P^\perp$, but not
  containing $P$;
\item a BLT-set of lines $\mathcal{L}$ of $\mathsf{W}(q)$.
\end{itemize}

For each element $\ell_i$ of $\mathcal{L}$, let $\pi_i$ be the plane
spanned by $\ell_i$ and $P$. Then we construct a generalised quadrangle as
follows:

\begin{table}[H]
\begin{center}
\begin{tabular}{l|l}
\hline
Points & Lines \\
\hline
-- points of $\mathsf{PG}(5,q)$ not in $P^\rho$&
-- totally iso. planes not contained in $P^\rho$ and \\
-- lines of $\mathsf{PG}(5,q)$ not incident&meeting some $\pi_i$ in a line\\
with $P$ but contained in some $\pi_i$&-- the planes $\pi_i$\\
-- the point $P$.&\\
\hline
\end{tabular}
\end{center}
\caption{The points and lines of the elation generalised quadrangle arising
from a BLT-set of lines of $\mathsf{W}(q)$. Incidence is inherited from
that of $\mathsf{PG}(5,q)$.}
\end{table}

Kantor \cite[Lemma]{Kantor91} showed that a Kantor family of the flock
elation group that is constructed from a $q$-clan, gives rise to a BLT-set
of lines of $\mathsf{W}(q)$. We show in Section \ref{proofs}, that for $q$
prime, any Kantor family of a flock elation group gives rise to a BLT-set
of lines of $\mathsf{W}(q)$, and the resulting flock quadrangle obtained by
the Knarr construction is isomorphic to the elation generalised quadrangle
arising from the given Kantor family.

\subsection*{Special and extraspecial groups}

A finite $p$-group $G$ is \textit{special}, if its centre, its derived
subgroup, and its Frattini subgroup coincide.  Moreover, we say that a
special group is \textit{extraspecial} if its centre is cyclic of prime
order.  The exponent of a special group is either $p$ or $p^2$. Further,
the order of an extraspecial group is of the form $p^{2m+1}$, where $m$ is
a positive integer. For each such $m$ there are, up to isomorphism,
precisely two extraspecial groups of order $p^{2m+1}$, one with exponent
$p$, and another one with exponent $p^2$~\cite[\S 8]{Aschbacher}.  The
elation groups of the flock quadrangles of order $(p^2,p)$ are extraspecial
of exponent $p$ (see \cite{Payne89}).

Here we recall a few facts about extraspecial $p$-groups which can be
readily found in \cite[\S 8]{Aschbacher}.  The quotient group $E/Z(E)$ is
an elementary abelian $p$-group forming a vector space $V$ over
$\GF(p)$. Moreover, the map from $V^2$ to $Z(E)$ defined by
$$\langle Z(E)x, Z(E)y\rangle = [x,y]$$
defines an alternating form on $V$. Thus, we obtain the generalised
quadrangle $\mathsf{W}(p)$ where the totally isotropic subspaces correspond
to abelian subgroups of $E$ properly containing $Z(E)$.  

%
%

\section{Kantor Families for $p$-Groups of Order $p^5$}\label{sectp5}

Recall that the elation group of a generalised quadrangle of order
$(p^2,p)$, $p$ prime, has order $p^5$. Thus we provide in this section some
powerful tools which will enable us to prove Theorem~\ref{Csaba}.

\begin{lem}\label{nicesubGQ}
  Let $(G,\mathcal{F},\mathcal{F}^*)$ be a Kantor family giving rise to an
  elation generalised quadrangle $\mathcal{S}$ of order $(s,t)$.  Suppose
  that $H$ is a subgroup of $G$ of order $t^3$ such that for all
  $A\in\mathcal{F}$ and $A^*\in\mathcal{F}^*$, we have
$$|A^* \cap H|\ge t^2\text{ and }|A \cap H|\ge t.$$ Then
$$(\{A\cap H:A\in\mathcal{F}\},\{A^*\cap H:A^*\in\mathcal{F}^*\})$$
is a Kantor family for $H$ giving rise to an elation generalised
quadrangle of order~$t$.
\end{lem}

\begin{proof}
  Suppose $A$ and $B$ are a pair of distinct elements of $\mathcal{F}$, and
  let $A^*$ and $B^*$ be the respective elements of $\mathcal{F}^*$ such
  that $A\le A^*$ and $B\le B^*$. Since $A$ and $B^*$ intersect trivially,
  we have that
$$|A^*H|\ge|A^*(B\cap H)|=\frac{|A^*||B\cap H|}{|A^*\cap B\cap H|} \ge
st^2.$$ Therefore $$|A^*\cap H|=|A^*||H|/|A^*H|\le t^2,$$ and so $A^*$ and
$H$ intersect in $t^2$ elements, for all $A^*\in\mathcal{F}^*$.  Similarly,
$$|AH|\ge|A(B^*\cap H)|=|A||B^*\cap H|\ge st^2$$ and so $|A \cap H|=t$, for all
$A\in\mathcal{F}$. The ``triple'' and ``tangency'' conditions follow from
those in $(G,\f,\ff)$.
\end{proof}

\begin{thm}\label{Brownetal}
  Let $p$ be an odd prime. A generalised quadrangle of order $(p^2,p)$ with
  an elation subquadrangle of order $p$ is isomorphic to $\mathsf{H}(3,p^2)$.
  Moreover, the subquadrangle here is isomorphic to $\mathsf{W}(p)$ and so
  is not a translation generalised quadrangle.
\end{thm}

\begin{proof}
  Let $\mathcal{S}$ be a generalised quadrangle of order $(p^2,p)$ with a
  subquadrangle $\mathcal{S}'$ of order $p$.  By \cite{BTvM}, a
  generalised quadrangle of order $p$ is either isomorphic to
  $\mathsf{W}(p)$ or $\mathsf{Q}(4,p)$. Now every line of our given
  generalised quadrangle of order $(p^2,p)$ induces a spread of the
  subquadrangle; but $\mathsf{Q}(4,p)$ has no spreads for $p$ odd (see
  \cite[3.4.1(i)]{FGQ}).  Therefore, $\mathcal{S}'$ is isomorphic to
  $\mathsf{W}(p)$. It was proved by Brown \cite{Brown02}, and independently
  by Brouns, Thas, and Van Maldeghem \cite{BrounsTvM}, that if a
  generalised quadrangle $\mathcal{S}$ of order $(q,q^2)$ has a
  subquadrangle $\mathcal{S}'$ isomorphic to $\mathsf{Q}(4,q)$, and if in
  $\mathcal{S}'$ each ovoid $\mathcal{O}_X$ consisting of all of the points
  collinear with a given point $X$ of $\mathcal{S}\backslash\mathcal{S}'$
  is an elliptic quadric, then $\mathcal{S}$ is isomorphic to
  $\mathsf{Q}(5,q)$.  By a result of Ball, Govaerts, and Storme \cite{BGS},
  if $p$ is a prime then every ovoid of $\mathsf{Q}(4,p)$ is an elliptic
  quadric. Therefore, by dualising, we have that $\mathcal{S}$ is
  isomorphic to $\mathsf{H}(3,p^2)$.
\end{proof}

The reason why we have pointed out that the subquadrangle is not a
translation generalised quadrangle will become apparent in Section
\ref{proofthcsaba}. We obtain the following consequence of Theorem
\ref{Brownetal}.

\begin{lem}\label{filter}
  Let $p$ be a prime and let $(G,\mathcal{F},\mathcal{F}^*)$ be a Kantor
  family giving rise to an elation generalised quadrangle $\mathcal{S}$ of
  order $(p^2,p)$. Suppose that $H$ is a subgroup of $G$ of order $p^3$
  with the property that, for all $A^*\in\mathcal{F}^*$, we have $|A^* \cap
  H|\ge p^2$. Then $\mathcal{S}$ is isomorphic to $\mathsf{H}(3,p^2)$.
\end{lem}

\begin{proof}
  Let $A\in\f$ and $A^*\in\ff$ such that $A\le A^*$.  The condition
  $|A^*\cap H|\ge p^2$ implies that $A^*H\neq G$. This gives $AH\neq G$,
  and so $|A\cap H|\ge p$.  Now it follows from Lemma \ref{nicesubGQ} that
  $H$ gives rise to an elation subquadrangle $\mathcal{S}'$ of
  order $p$.  The remainder follows from Theorem \ref{Brownetal}.
\end{proof}

For $p$ odd, $\mathsf{W}(p)$ is not a translation generalised quadrangle,
which implies in the previous lemma that $H$ is non-abelian. The next
result gives more information about Kantor families for groups of order
$p^5$.

\begin{lem}\label{l1}
  Suppose that $G$ is a group with order $p^5$ and let $(\f,\ff)$ be a
  Kantor family of order $(p^2,p)$ for $G$. Then the following hold:
\begin{enumerate}
\item[(i)] None of the members of $\f$ is normal in $G$. In particular $G$
  is non-abelian.
\item[(ii)] If $G$ is not extraspecial and $H$ is a subgroup of $G$ of
  order $p^3$, then there is a subgroup $U$ of $G$ such that $|U|=p^3$ and
  $HU=G$.
\item[(iii)] $G$ is not generated by two elements.
\item[(iv)] The nilpotency class of $G$ is two.
\item[(v)] $G'$ is elementary abelian.
\end{enumerate}
\end{lem}
\begin{proof}
  If $G$ is an extraspecial group with order $p^5$, then properties (i),
  (iii), (iv), and (v) are valid for $G$, and so we may assume, for the
  entire proof, that $G$ is not extraspecial.

  (i) Assume by contradiction that $A\in \f$ is normal, and choose distinct
  $B,\ C\in\f\setminus\{A\}$. Then $AB$ is a subgroup of $G$ with order
  $p^4$ and so $AB\cap C=1$ is impossible, violating the triple condition.

  (ii) Let $H$ be a subgroup of $G$ with order $p^3$. Since the elation
  group of $\mathsf{H}(3,p^2)$ is extraspecial with exponent $p$,
  Lemma~\ref{filter} implies that there is $A^*\in\ff$ such that $|H\cap
  A^*|=p$, and so $HA^*=G$.

  (iii) Since $G/\Phi(G)$ is not cyclic, $|\Phi(G)|\le p^3$.  Further,
  $\Phi(G)U=G$ implies that $U=G$, and hence it follows from part~(ii) that
  $\Phi(G)\neq p^3$. Therefore we obtain that $|\Phi(G)|\le p^2$, and so a
  minimal generating set of $G$ has at least three elements.

  (iv) A group of order $p^5$ has nilpotency class at most $4$.  If the
  nilpotency class of $G$ is $4$, then $|G'|=|\Phi(G)|=p^3$, which is a
  contradiction by the previous paragraph.  We claim that the nilpotency
  class of $G$ is not three. Suppose by contradiction that it is three. In
  this case, as $G$ is not generated by 2 elements, $G/G'\cong C_p\times
  C_p\times C_p$ and $|G':\lcs G3|=|\lcs G3|=p$.  Choose $a,\ b\in G$ such
  that $\left<[a,b]\lcs G3\right>=G'/\lcs G3$. Let $c_1\in G$ such that
  $\left<aG',bG',c_1G'\right>=G/G'$.  Then there are $\alpha,\ \beta\in
  \F_p$ such that $[a,c_1]\equiv [a,b]^\alpha\pmod{\lcs G3}$ and
  $[b,c_1]\equiv [a,b]^\beta\pmod{\lcs G3}$. Set $c=c_1a^\beta
  b^{-\alpha}$.  Then $\left<aG',bG',cG'\right>=G/G'$ and
  $[a,c]\equiv[b,c]\equiv 1\pmod{\lcs G3}$; that is $[a,c],\ [b,c]\in\lcs
  G3$.  By the Hall-Witt identity, $[a,b,c]=[c,b,a][a,c,b]=1$. As $\lcs
  G3=\left<[a,b,a],[a,b,b],[a,b,c]\right>$, this implies that either
  $[a,b,a]\neq 1$ or $[a,b,b]\neq 1$. Hence the subgroup $\left<a,b\right>$
  has nilpotency class $3$ and order $p^4$ (see
  also~\cite[Corollary~2.2(i)]{sch}).

  Let $H=\left<c,G'\right>$. Clearly, $|H|=p^3$ and
  $G/H=\left<aH,bH\right>$.  Let $U$ be a subgroup of $G$ such that $HU=G$,
  and so $HU/H=G/H=\left<aH,bH\right>$. This shows that there are $h_1,\
  h_2\in H$ such that $ah_1,\ bh_2\in U$. Since $[a,h_1],\ [a,h_2]\in\lcs
  G3$ and $[a,b,h_1]=[a,b,h_2]=1$, we obtain that $[ah_1,bh_2]\lcs
  G3=G'/\lcs G3$ and either $[ah_1,bh_2,ah_1]\neq 1$ or
  $[ah_1,bh_2,bh_2]\neq 1$. Thus $U$ contains $G'$ and $U$ is a group of
  order at least $p^4$. This, however, is a contradiction, by part~(ii).
  Therefore the nilpotency class of $G$ is not three. Since, by part~(i),
  the nilpotency class of $G$ is not one, we obtain that the class of $G$
  must be two.

  (v) By~(iv), we only need to show that the exponent of $G'$ is $p$.
  By~\cite[5.2.5]{robinson}, the quotient $G'/\lcs G3=G'$, as an abelian
  group, is an epimorphic image of the tensor product $(G/G')\otimes_{\Z}
  (G/G')$, which implies that the exponent of $G'/\lcs G3=G'$ divides the
  exponent of $G/G'$.  As $G$ is not generated by two elements, the size,
  and hence the exponent, of $G'$ is at most $p^2$. However, if this
  exponent is $p^2$ then $G/G'\cong (C_p)^3$, which is impossible.
\end{proof}

The next lemma describes the case when either $G'$ or $\Phi(G)$ is small.

\begin{lem}\label{l2}
  Suppose that $G$ is a group with order $p^5$ and let $(\f,\ff)$ be a
  Kantor family for $G$.
\begin{enumerate}
\item[(i)] If $|G'|=p$ then all members of $\f\cup \ff$ are abelian.
\item[(ii)] If $|\Phi(G)|=p$ then all members of $\f\cup \ff$ are
  elementary abelian. Moreover, if $p$ is odd, then in this case, $G$ has
  exponent $p$.
\item[(iii)] If $p$ is odd and $G$ is extraspecial, then $G$ has exponent
  $p$ and all members of $\f\cup\ff$ are elementary abelian.
\end{enumerate}
\end{lem}

\begin{proof}
  (i) Let us first assume that $|G'|=p$.  It suffices to prove, for all
  $A^*\in\ff$, that $A^*$ is abelian. We argue by contradiction and assume
  that $A^*\in\ff$ is not abelian.  In this case the derived subgroup
  $(A^*)'$ of $A^*$ is non-trivial, and, as $(A^*)'\le G'$, we obtain that
  $(A^*)'=G'$. Let $A\in\f$ such that $A\le A^*$. Then $A$ is a maximal
  subgroup of $A^*$, and so $(A^*)'=G'\le A$. Thus $A$ is normal in $G$,
  which is impossible by Lemma~\ref{l1}(i).  Therefore $A^*$ is abelian, as
  claimed.

  (ii) The assertion that the members of the Kantor family are elementary
  abelian can be proved by substituting $\Phi(A^*)$ in the place of
  $(A^*)'$ and $\Phi(G)$ in the place of $G'$ in the previous paragraph.
  Let $p$ be an odd prime. In this case, as $|G'|=p$, the elements of $G$
  with order $p$ form a subgroup $\Omega(G)$ of $G$. Let $A\in \f$ and
  $B^*\in\ff$ such that $A\cap B^*=1$. In this case $AB^*=G$ and $A,\
  B^*\le\Omega(G)$. Therefore $G=\Omega(G)$, which amounts to saying that
  $G$ has exponent $p$.

  Part~(iii) follows immediately from part~(ii).
\end{proof}

The following lemma is a generalisation of \cite[Lemma]{Kantor91}.

\begin{lem}\label{extraspecial}
  Let $p$ be an odd prime and let $(\mathcal{F},\mathcal{F}^*)$ be a
  Kantor family for an extraspecial group $E$ of order $p^5$.  Then the
  image of $\mathcal{F}^*$ in $E/Z(E)$ corresponds to a BLT-set of lines of
  $\mathsf{W}(p)$.
\end{lem}

\begin{proof}
  First note that by Lemma \ref{l2}(iii), all the members of $\ff$ are abelian
  and hence each $A^*\in\ff$ induces an abelian subgroup of $E/Z(E)$, and so
a totally isotropic line of the associated $\mathsf{W}(p)$ geometry.
Therefore, every member of $\ff$ contains $Z(E)$.
  Suppose by way of contradiction that there is a line of $\mathsf{W}(p)$
  concurrent with three elements of $\mathcal{L}=
  \{A^*/Z(E):A^*\in\mathcal{F}^*\}$. Then there exists an abelian subgroup
  $H$ of $E$ of order $p^3$, and three elements $A^*,\ B^*,\ C^*$ of
  $\mathcal{F}^*$ such that $H$ intersects each of these elements in a
  subgroup of order $p^2$ properly containing $Z(E)$ (note: $H$ contains
  $Z(E)$). Let $A$, $B$, $C$ be the unique elements of $\mathcal{F}$
  contained in $A^*$, $B^*$, $C^*$ respectively. Now $(H\cap B)Z(E)$ is
  contained in $B^*$ and so $A\cap (H\cap B)Z(E)=1$.  Also, we have that
  $|H\cap A|= p$ as $p^2=|H\cap A^*|=|(H\cap A)Z(E)|=|H\cap A||Z(E)|$
  (similarly, $|H\cap B|= p$). Thus
\begin{align*}
|(H\cap A)(H\cap B)Z(E)|&=\frac{|H\cap A||(H\cap B)Z(E)|}{|H\cap A\cap (H\cap B)Z(E)|}\\
&=|H\cap A||(H\cap B)Z(E)|\\
&=|H\cap A||H\cap B||Z(E)|\\
&\ge p^3
\end{align*}
and so one can see that $H=(H\cap A)(H\cap B)Z(E)$. So
$$C^*\cap H=(C\cap (H\cap A)(H\cap B))Z(E)$$
and by the condition $AB\cap C=1$, we have that $C^*\cap H=Z(E)$, giving us
the desired contradiction. Therefore, $\mathcal{L}$ is a BLT-set of lines
of $\mathsf{W}(p)$.
\end{proof}

%
%

\section{The Proof of Theorem \ref{Csaba}}\label{proofthcsaba}

In this section we prove Theorem~\ref{Csaba}.  
By Lemma~\ref{l2}(iii), an extraspecial
group with order $p^5$ and exponent $p^2$ does not admit a Kantor family
with order $(p^2,p)$. Hence we may assume, for a proof by contradiction, that:

\begin{quote} 
{\em $G$ is a group of order $p^5$ 
 and $(\f,\ff)$ is a Kantor family for $G$ with order $(p^2,p)$.}
\end{quote} 
Our aim is to derive a contradiction.  First note that Lemma~\ref{l1}
implies that one of the following must hold:
\begin{enumerate}
\item[(I)] $G/G'\cong C_p\times C_p\times C_{p^2}$ and $G'\cong C_p$;
\item[(II)] $G/G'\cong (C_p)^3$ and $G'\cong (C_p)^2$;
\item[(III)] $G/G'\cong (C_p)^4$ and $G'\cong C_p$.
\end{enumerate}
We show, case by case, that none of the above possibilities can occur.  We
let $Z$ denote the centre of $G$.

\medskip\noindent Case~(I).  Using the argument in the proof of
Lemma~\ref{l1}(iv), we can choose generators $a,\ b,\ c$ of $G$ such that
$G'=\left<[a,b]\right>$ and $c\in Z$. It also follows that
$Z=\left<z,\Phi(G)\right>$, and so $|Z|=p^3$. By Lemma~\ref{l2}(ii), all
members of $\ff$ must be abelian and so~\cite[Theorem~3.2 and
Lemma~2.2]{Hach96} imply that the subgroups $A\cap Z$ and $A^*\cap Z$ with
$A\in \f$ and $A^*\in\ff$ form a Kantor family for $Z$ with order
$p$. This, however, contradicts Theorem \ref{Brownetal}, since the
subquadrangle here is not a translation generalised quadrangle (n.b., $Z$
is abelian). Hence case~(I) cannot occur.

\medskip\noindent Case~(II).  First we claim that it is possible to choose
the generators $x$, $y$, and $z$ of $G$ such that
$G'=\left<[x,y],[x,z]\right>$ and $[y,z]=1$.  Let $x,\ y,\ z$ be generators
of $G$. Then $G'=\left<[x,y],[x,z],[y,z]\right>$.  Since $G'\cong (C_p)^2$
we have that there are $\alpha,\ \beta,\ \gamma\in\F_p$ such that at least
one of $\alpha,\ \beta,\ \gamma$ is non-zero and
$[x,y]^\alpha[x,z]^\beta[y,z]^\gamma=1$. If $\alpha=\beta=0$ then
$\gamma\neq 0$, and $[y,z]=1$ follows. If $\alpha=0$ and $\beta\neq 0$ then
$[x^\beta y^\gamma,z]=1$. Now replacing $x$ by $x^\beta y^\gamma$ we find
that in the new generating set $[x,z]=1$ holds. Similarly, if $\alpha\neq
0$ and $\beta=0$ then $[y,x^{-\alpha}z^\gamma]=1$ and replacing $x$ by
$x^{-\alpha}z^\gamma$ we obtain that $[x,y]=1$ holds in the new generating
set. Finally if $\alpha\beta\neq 0$, then we replace $x$ by
$x^{\beta/\alpha}y^{\gamma/\alpha}$ and $y$ by $yz^{\beta/\alpha}$ to
obtain that $[x,y]=1$. Thus, after applying one of the substitutions above
and possibly renaming the generators, $[y,z]=1$ holds, and the claim is
valid.

We continue by verifying the following claim: if $H$ is a subgroup in $G$
with order $p^2$ and $H\cap Z=1$ then there are $c,\ d\in Z$ such that
$H=\left<yc,zd\right>$.

Assume that $H$ is a subgroup or order $p^2$ that does not intersect
$Z$. Then $HZ/Z\cong H/(H\cap Z)=H$ and so $H\cong C_p\times C_p$.  In
particular $H$ can be generated by two elements of the form $u=x^{\alpha_1}
y^{\beta_1} z^{\gamma_1} c_1$ and $v=x^{\alpha_2} y^{\beta_2} z^{\gamma_2}
c_2$ where $\alpha_i,\ \beta_i,\ \gamma_i\in\F_p$, $c_i\in Z$ and
$\left<uZ,vZ\right>\cong C_p\times C_p$. Since $[u,v]=1$ we obtain that
$$
1=[u,v]=
[x^{\alpha_1} y^{\beta_1} z^{\gamma_1} c_1,x^{\alpha_2} y^{\beta_2}
z^{\gamma_2}
c_2]=[x,y]^{\alpha_1\beta_2-\alpha_2\beta_1}[x,z]^{\alpha_1\gamma_2-\alpha_2\gamma_1}. 
$$
Thus $\alpha_1\beta_2-\alpha_2\beta_1=\alpha_1\gamma_2-\alpha_2\gamma_1=0$.
Note that these two expressions can be viewed as determinants of suitable
$(2\times 2)$-matrices.  If $(\alpha_1,\alpha_2)\neq (0,0)$ then the
vectors $(\beta_1,\beta_2)$ and $(\gamma_1,\gamma_2)$ are both multiples of
$(\alpha_1,\alpha_2)$ and so the matrix
$$
\left(
\begin{array}{cc}
\alpha_1 &\alpha_2\\
\beta_1 &\beta_2\\
\gamma_1 & \gamma_2
\end{array}\right)
$$
has row-rank~1. Since the row-rank of a matrix is the same as the
column-rank, this also shows that the vector $(\alpha_2,\beta_2,\gamma_2)$
is a multiple of the vector $(\alpha_1,\beta_1,\gamma_1)$ and so $uZ=vZ$
which gives $HZ/Z\cong C_p$: a contradiction. Thus
$(\alpha_1,\alpha_2)=(0,0)$; thus $u=y^{\beta_1}z^{\gamma_1}c_1$ and
$v=y^{\beta_2}z^{\gamma_2}c_2$. Since $\left<uZ,vZ\right>\cong C_p\times
C_p$ we must have that $\beta_1\gamma_2-\beta_2\gamma_1\neq 0$. Also, if
$\beta_1,\ \beta_2=0$ then $HZ/Z\cong C_p$ and so we may assume that
$\beta_1\neq 0$.  Change $v$ to $u^{-\beta_2/\beta_1}v$; then
$\left<u,v\right>=H$ and $v$ is of the form $z^{\gamma}d'$ where $d'\in Z$.
Now change $u$ to $uv^{-\gamma_1/\gamma_2}$.  Then $\left<u,v\right>=H$
still holds and now $u$ is of the form $y^{\beta}c'$ where $c'\in Z$.  Now
$u^{\beta^{-1}}$ and $v^{\gamma^{-1}}$ are as required.

Let us now prove that $G$ does not admit a Kantor family.  We argue by
contradiction and assume that $(\f,\ff)$ is a Kantor family of order
$(p^2,p)$ for $G$.  If $A,\ B$ are distinct elements of $\f$ such that
$A\cap Z=B\cap Z=1$ then the claim above implies that $[A,B]=1$, and so
$AB$ is a subgroup of $G$ with order $p^4$. Thus, if
$C\in\f\setminus\{A,B\}$, then $AB\cap C\neq 1$, which contradicts the
triple condition. Thus $\f$ has at most one member that avoids the
centre. Let us suppose now $A,\ B,\ C$ are pairwise distinct members of
$\f$ such that $A\cap Z$, $B\cap Z$, and $C\cap Z$ are non-trivial. As
$A\cap B=A\cap C=B\cap C=1$, we obtain that $|A\cap Z|=|B\cap Z|=|C\cap
Z|=p$ and that $A\cap Z$, $B\cap Z$, $C\cap Z$ are three distinct subgroups
of $Z$.  This, however implies that $Z=(A\cap Z)(B\cap Z)$, and, in turn,
that $C\cap Z\le (A\cap Z)(B\cap Z)$, which violates the triple condition.

The argument in the last paragraph implies that at most two members of $\f$
can intersect $Z$ non-trivially, and at most one member of $\f$ can avoid
the centre. Thus $|\f|\le 3$, which is a contradiction as $p$ is odd and
$|\f|=p+1$. Therefore case~(II) is impossible.

\medskip\noindent Case~(III).  As $G$ is not extraspecial, $|Z|=p^3$, and
Lemma~\ref{l1}(iv) implies that the members of $\ff$ are abelian.  In this
case~\cite[Theorem~3.2, Lemmas~2.1 and~2.2]{Hach96} show that the subgroups
$A\cap Z$ and $A^*\cap Z$ (with $A\in\f$ and $A^*\in\ff$) form a Kantor
family of order $p$ for $Z$. However, we have a contradiction to Theorem
\ref{Brownetal} since the associated subquadrangle of order $p$ is not a
translation generalised quadrangle.

As none of the possibilities listed at the beginning of the section can
occur, Theorem~\ref{Csaba} must hold.

%
%

\section{The Proof of Theorem \ref{Theorem1}}\label{proofs}

Here we prove Theorem \ref{Theorem1}, but first we show that applying 
the Knarr construction to a
BLT-set of lines arising from a Kantor family $(\mathcal{F},\mathcal{F}^*)$
of the flock elation group results in an elation generalised quadrangle
isomorphic to that directly associated to $(\mathcal{F},\mathcal{F}^*)$.

\begin{thm}\label{ShultThasKnarr}
  Let $G$ be the flock elation group of order $p^5$, $p$ odd, and suppose
  that $G$ admits a Kantor family $(\mathcal{F},\mathcal{F}^*)$ giving rise
  to an elation generalised quadrangle $\mathcal{E}$.  Consider the BLT-set
  of lines $\mathcal{L}$ of $\mathsf{W}(p)$ obtained by taking the image of
  $\mathcal{F}^*$ under the natural projection map from $G$ onto $G/Z(G)$.
  Then the flock quadrangle arising from $\mathcal{L}$ via the Knarr
  construction is equivalent to $\mathcal{E}$.
\end{thm}

\begin{proof}
First note that $G$ is extraspecial of exponent $p$, and 
observe that the matrices of the form
$$\left(\begin{smallmatrix}
    1&a&b&c&d&e\\
    0&1&0&0&0&d\\
    0&0&1&0&0&c\\
    0&0&0&1&0&-b\\
    0&0&0&0&1&-a\\
    0&0&0&0&0&1
\end{smallmatrix}\right), \quad a,b,c,d,e\in\GF(p)$$
define a representation of $G$ into the symplectic group
$\mathrm{PSp}(6,p)$ with its associated null polarity given by the matrix
$$\left(\begin{smallmatrix}
0&0&0&0&0&1\\
0&0&0&0&1&0\\
0&0&0&1&0&0\\
0&0&-1&0&0&0\\
0&-1&0&0&0&0\\
-1&0&0&0&0&0
\end{smallmatrix}\right).$$

Moreover, the centre of $G$ consists only of those upper triangular
matrices with zeros everywhere above the diagonal except possibly the top
right corner, and $G$ fixes the projective point $P$ represented by
$(1,0,0,0,0,0)$. Hence $G$ induces an action on the quotient
$P^\perp/P\equiv \mathsf{W}(p)$.  It is not difficult to show that the
right coset action of $G$ on $G/Z(G)$ is permutationally isomorphic to the
action of $G$ on $P^\perp/P$ (as a projective right-module). To be more
specific, the representatives of $G/Z(G)$ are in a bijection with matrices
of the form
$$\left(\begin{smallmatrix}
1&a&b&c&d&0\\
0&1&0&0&0&d\\
0&0&1&0&0&c\\
0&0&0&1&0&-b\\
0&0&0&0&1&-a\\
0&0&0&0&0&1
\end{smallmatrix}\right), \quad a,b,c,d\in\GF(p)$$
and $P^\perp/P$ can naturally be identified with vectors of the form
$(0,a,b,c,d,1)$. Thus we have a bijection from $G/Z(G)$ onto $P^\perp/P$
given by
$$\left(\begin{smallmatrix}
1&a&b&c&d&0\\
0&1&0&0&0&d\\
0&0&1&0&0&c\\
0&0&0&1&0&-b\\
0&0&0&0&1&-a\\
0&0&0&0&0&1
\end{smallmatrix}\right)\mapsto P+(0,a,b,c,d,1)$$
such that the right coset action of $G$ is equivalent to the right-module
action of $G$ on $P^\perp/P$.

Let $(\mathcal{F},\mathcal{F^*})$ be a Kantor family for $G$ and let
$\mathcal{E}$ be the associated elation generalised quadrangle with points
\begin{enumerate}
\item[(i)] elements of $g$,
\item[(ii)] right cosets $A_i^*g$ of elements of $\mathcal{F}^*$,
\item[(iii)] $\infty$,
\end{enumerate}
and lines
\begin{enumerate}
\item[(a)] right cosets $A_ig$ of elements of $\mathcal{F}$,
\item[(b)] symbols $[A_i]$ where $A_i\in\mathcal{F}$.
\end{enumerate}

Let $Q=(0,0,0,0,0,1)$ and note that $Q$ is opposite to $P$.  Let
$\mathcal{K}$ be the flock quadrangle associated to $\mathcal{L}$
constructed from the point $P$, and define a map from $\mathcal{E}$ to
$\mathcal{K}$ as follows:
$$
\infty\mapsto P,\quad
[A_i]\mapsto \pi_i,\quad
A_i^*g\mapsto z_i^g,\quad
A_ig\mapsto M_i^g,\quad
g\mapsto Q^g.$$

We will show that this map defines an isomorphism of generalised
quadrangles. Since the action of $G$ on $P^\perp/P$ is permutationally
isomorphic to the right coset action of $G$ on $G/Z(G)$, we have that the
stabiliser of the subspace corresponding to a subgroup $H$ containing
$Z(G)$ is just $H$ itself. Therefore $A_i$ fixes $z_i$ and $A_i^*$ fixes
$M_i$ (for all $i$), and so the map above is well-defined. Now we verify
that the four types of incidences are compatible:

\medskip\noindent\textbf{Incidence of $\mathbf{\infty}$ and $\mathbf{[A_i]}$:} 
It is clear that $P\sim \pi_i$ for all $i$.

\medskip\noindent\textbf{Incidence of $\mathbf{A_i^*g}$ and
  $\mathbf{[A_i]}$:} We want to show that $\pi_i\sim z_i^g$ given we know
that $\pi_i\sim z_i$. Now $G$ fixes every subspace of $P^\perp$ on $P$, and
hence $G$ fixes $\pi_i$.  Therefore $z_i^g\sim\pi_i^g=\pi$ (n.b., $g$ is a
collineation).

\medskip\noindent\textbf{Incidence of $\mathbf{A_i^*g}$ and
  $\mathbf{A_ih}$:} So $A_ih\subset A_i^*g$. We want to show that
$M_i^h\sim z_i^g$. By definition, $z_i$ is the unique line of $\pi_i$ (not
on $P$) which is on a plane $M_i$ on $Q$. We know that $M_i\sim z_i$. Since
$A_ih\subset A_iZg$, then there exists an element $e$ of $Z(G)$ such that
$hg^{-1}e\in A_i$. It suffices to show that $M_i^{hg^{-1}}\sim z_i$.  Now
$A_i$ fixes $M_i$ and so $M_i^{hg^{-1}}=M_i^{e^{-1}}$.  Now $e^{-1}$ fixes
$z_i$ and so $M_i^{hg^{-1}}\sim z_i$.

\medskip\noindent\textbf{Incidence of $\mathbf{g}$ and $\mathbf{A_ig}$:} It
is clear that $G$ acts regularly on the points opposite $P$. Since for all
$i$ we have $Q\sim M_i$, it follows that $Q^g\sim M_i^g$.

\medskip

Therefore, the flock quadrangle arising from $\mathcal{L}$ via the Knarr
construction is equivalent to $\mathcal{E}$.
\end{proof}

\subsection*{Theorem \ref{Theorem1} and its proof}

\begin{quote}
  {\em An elation generalised quadrangle of order $(s,p)$, with $p$
    prime, is a flock quadrangle, isomorphic to $\mathsf{Q}(4,p)$,
     or isomorphic to $\mathsf{W}(p)$.}
\end{quote}

\begin{proof}
  Let $\mathcal{S}$ be an elation generalised quadrangle of order $(s,p)$,
  where $p$ is prime, and suppose that $(G,\mathcal{F},\mathcal{F}^*)$ is
  the corresponding Kantor family. By \cite{BTvM}, we may assume
  $s=p^2$, and so $G$ has order $p^5$.
  By Theorem \ref{Csaba}, $G$ must be extraspecial. Now
  the Frattini subgroup of $G$ has order $p$ and so has nontrivial
  intersection with every subgroup of $G$ that has order at least $p^3$.
  Hence $Z(G)$ is contained in every element of $\mathcal{F}^*$.
  Therefore, by Lemma \ref{extraspecial} and Theorem \ref{ShultThasKnarr},
  our generalised quadrangle $\mathcal{S}$ is a flock quadrangle.
\end{proof}

\section*{Acknowledgment}
The authors would like to thank Stan Payne for his feedback on this
work. This research formed part of an Australian Research Council Discovery
Grant project that was undertaken at the University of Western
Australia. The first author was supported by a Marie Curie Incoming
International Fellowship within the 6th European Community Framework
Programme, contract number: MIIF1-CT-2006-040360.  The third author was
supported by the Hungarian Scientific Research Fund (OTKA) grant F049040;
he is also grateful to Hendrik Van Maldeghem for funding his visit to the
workshop {\em Groups and Buildings 2007} held in Ghent.

\end{document}